\documentclass[12pt, twoside, eqno]{article}
\usepackage{latexsym}
\usepackage{amssymb}
\usepackage{amsfonts}
\begin{document}

\noindent {\bf \large A characterization of regular, intra-regular, 
left quasi-regular and semisimple hypersemigroups in terms of fuzzy 
sets}
\bigskip

\noindent{\bf Niovi Kehayopulu}\\
June 16, 2016\bigskip{\small

\noindent{\bf Abstract.} We prove that an hypersemigroup $H$ is 
regular if and only, for any fuzzy subset $f$ of $H$, we have 
$f\preceq f\circ 1\circ f$ and it is intra-regular if and only if, 
for any fuzzy subset $f$ of $H$, we have $f\preceq 1\circ f\circ 
f\circ 1$. An hypersemigroup $H$ is left (resp. right) quasi-regular 
if and only if, for any fuzzy subset $f$ of $S$ we have $f\preceq 
1\circ f\circ 1\circ f$ (resp. $f\preceq f\circ 1\circ f\circ 1)$ and 
it is semisimple if and only if, for any fuzzy subset $f$ of $S$ we 
have $f\preceq 1\circ f\circ 1\circ f\circ 1$. The characterization 
of regular and intra-regular hypersemigroups in terms of fuzzy 
subsets are very useful for applications.\smallskip

semisimple}
\section{Introduction}In our paper in [2] we gave, among others, a 
characterization of regular ordered semigroups in terms of fuzzy 
subsets which is very useful for applications. Using that equivalent 
definition of regular semigroups many known results on semigroup 
(without order) or on ordered semigroups can be drastically 
simplified. In [1], we characterized the left quasi-regular and 
semisimple ordered semigroups in terms of fuzzy sets. In the present 
paper we characterize the regular, intra-regular, left quasi-regular 
and semisimple hypersemigroups using fuzzy sets. According to the 
equivalent definition of regularity and intra-regularity given in the 
present paper many results on hypersemigroups can be drastically 
simplified. The paper has been inspired by our paper in [2], and the 
aim is to show the way we pass from fuzzy semigroups to fuzzy 
hypersemigroups. In fact, the results on semigroups or ordered 
semigroups can be transferred to hypersemigroups in the way indicated 
in the present paper.
\section{Main results} An {\it hypergroupoid} is a nonempty set $H$ 
with an hyperoperation $$\circ : H\times H \rightarrow {\cal P}^*(H) 
\mid (a,b) \rightarrow a\circ b$$on $H$ and an operation $$* : {\cal 
P}^*(H)\times {\cal P}^*(H) \rightarrow {\cal P}^*(H) \mid (A,B) 
\rightarrow A*B$$ on ${\cal P}^*(H)$ (induced by the operation of 
$H$) such that $$A*B=\bigcup\limits_{(a,b) \in\,A\times B} {(a\circ 
b)}$$ for every $A,B\in {\cal P}^*(H)$ (${\cal P}^*(H)$ denotes the 
set of nonempty subsets of $H$).

The operation ``$*$" is well defined. Indeed: If $(A,B)\in {\cal 
P}^*(H) \times {\cal P}^*(H)$, then $A*B=\bigcup\limits_{(a,b) 
\in\,A\times B} {(a\circ b)}$. For every $(a,b)\in A\times B$, we 
have $(a,b)\in H\times H$, then $(a\circ b)\in {\cal P}^*(H)$, thus 
we get $A*B\in {\cal P}^*(H)$. If $(A,B),(C,D)\in {\cal P}^*(H)\times 
{\cal P}^*(H)$ such that $(A,B)=(C,D)$, then 
$$A*B=\bigcup\limits_{(a,b) \in\,A\times B} {(a\circ 
b)}=\bigcup\limits_{(a,b) \in\,C\times D} {(a\circ b)}=C*D.$$
As the operation ``$*$" depends on the hyperoperation ``$\circ$", an 
hypergroupoid can be also denoted by $(H,\circ)$ (instead of 
$(H,\circ,*)$).

If $H$ is an hypergroupoid then, for any $x,y\in H$, we have $x\circ 
y=\{x\}*\{y\}$. Indeed,$$\{x\}*\{y\}=\bigcup\limits_{u \in \{ x\} ,v 
\in \{ y\} } {u \circ v = x \circ y}.$$

An hypergroupoid $H$ is called {\it hypersemigroup} if $$(x\circ 
y)*\{z\}=\{x\}*(y\circ z)$$ for every $x,y,z\in H$. Since $x\circ 
y=\{x\}*\{y\}$ for any $x,y\in H$, an hypergroupoid $H$ is an 
hypersemigroup if and only if, for any $x,y,z\in H$, we 
have$${\Big(}\{x\}*\{y\}{\Big)}*\{z\}=\{x\}*{\Big(}\{y\}*\{z\}{\Big)}.$$

Following Zadeh, if $(H,\circ)$ is an hypergroupoid, we say that $f$ 
is a fuzzy subset of $H$ (or a fuzzy set in $H$) if $f$ is a mapping 
of $H$ into the real closed interval $[0,1]$ of real numbers, that is 
$f : H \rightarrow [0,1]$. For an element $a$ of $H$, we denote by 
$A_a$ the subset of $H\times H$ defined as follows:$$A_a:=\{(y,z)\in 
H\times H \mid a\in y\circ z\}.$$For two fuzzy subsets $f$ and $g$ of 
$H$, we denote by $f\circ g$ the fuzzy subset of $H$ defined as 
follows:
$$f \circ g: H \to [0,1]\,\,\,a \to \left\{ \begin{array}{l}
\bigvee\limits_{(y,z) \in {A_a}} {\min \{ f(y),g(z)\} 
\,\,\,\,if\,\,\,{A_a} \ne \emptyset } \\
\,\,\,\,0\,\,\,\,\,\,if\,\,\,\,{A_a} = \emptyset. \,\,\,\,\,\,\,\,\,
\end{array} \right.$$Denote by $F(H)$ the set of all fuzzy subsets of 
$H$ and by ``$\preceq$" the order relation on $F(H)$ defined 
by:$$f\preceq g \;\Longleftrightarrow\; f(x)\le g(x) \mbox { for 
every } x\in H.$$We finally show by 1 the fuzzy subset of $H$ defined 
by:$$1: H \rightarrow [0,1] \mid x \rightarrow 1(x):=1.$$Clearly, the 
fuzzy subset 1 is the greatest element of the ordered set 
$(F(H),\preceq)$ (that is, $1\succeq f$ $\forall f\in F(H)$).

We denote the hyperoperation on $H$ and the multiplication between 
the two fuzzy subsets of $H$ by the same symbol (no confusion is 
possible).\medskip

For an hypergroupoid $H$, we denote by $f_a$ the fuzzy subset of $H$ 
defined by:$${f_a}:H \to [0,1] \mid x \to {f_a}(x): = \left\{ 
\begin{array}{l}
1\,\,\,\,\,if\,\,\,\,x = a\\
0\,\,\,\,if\,\,\,\,x \ne a.
\end{array} \right.$$
The following proposition, though clear, plays an essential role in 
the theory of hypergroupoids.\medskip

\noindent{\bf Proposition 1.} {\it Let $(H,\circ)$ be an 
hypergroupoid, $x\in H$ and $A,B\in {\cal P}^*(H)$. Then we have the 
following:

$(1)$ $x\in A*B$ $\Longleftrightarrow$ $x\in a\circ b$ for some $a\in 
A$, $b\in B$.

$(2)$ If $a\in A$ and $b\in B$, then $a\circ b\subseteq A*B$.} 
\medskip

\noindent{\bf Lemma 2}. {\it If $H$ is an hypergroupoid and $A,B,C$ 
nonempty subsets of $H$, then $A\subseteq B$ implies $A*C\subseteq 
B*C$ and $C*A\subseteq C*B$}.\medskip

\noindent{\bf Lemma 3.} (cf. also [3; Proposition 9]) {\it If $H$ is 
an hypersemigroup, then the set of all fuzzy subsets of $H$ is a 
semigroup.}\\According to this lemma, for any fuzzy subsets $f,g,h$ 
of $H$, we write

$(f\circ g)\circ h=f\circ (g\circ h):=f\circ g\circ h$.\medskip

\noindent{\bf Definition 4.} An hypersemigroup $H$ is called {\it 
regular} if for every $a\in H$ there exists $x\in H$ such that$$a\in 
(a\circ x)*\{a\}.$$Equivalent Definitions:

(1) $a\in \{a\}*H*\{a\}$ for every $a\in H$.

(2) $A\subseteq A*H*A$ for every $A\subseteq H$.\medskip

\noindent{\bf Theorem 5.} {\it An hypersemigroup H is regular if and 
only if, for any fuzzy subset f of H, we have$$f\preceq f\circ 1\circ 
f.$$}{\bf Proof.} $\Longrightarrow$. Let $f$ be a fuzzy subset of $H$ 
and $a\in H$. Since $H$ is regular, there exists $x\in H$ such that 
$a\in (a\circ x)*\{a\}.$ Then, By Proposition 1(1), there exists 
$u\in a\circ x$ such that $a\in u\circ a$. Since $(u,a)\in A_a$, we 
have $A_a\not=\emptyset$ and$$(f \circ 1 \circ f)(a): = 
\bigvee\limits_{(y,z) \in {A_a}} {\min } \{ (f \circ 1)(y),f(z)\}  
\ge \min \{ (f \circ 1)(u),f(a)\}.$$Since $(a,x)\in A_u$, we have 
$A_u\not=\emptyset$ and$$(f \circ 1)(u): = \bigvee\limits_{(y,z) \in 
{A_u}} {\min \{ f(y),1(z)\}  \ge \min \{ f(a),1(x)\}  = f(a).}$$Thus 
we have$$(f\circ 1\circ f)(a)\ge \min\{(f\circ 1)(u),f(a)\}\ge 
\min\{f(a),f(a)\}=f(a).$$$\Longleftarrow$. Let $a\in H$. Since $f_a$ 
is a fuzzy subset of $H$, by hypothesis, we have $1=f_a(a)\le 
(f_a\circ 1\circ f_a)(a)$. Since $f_a\circ 1\circ f_a$ is a fuzzy 
subset of $H$, we have $(f_a\circ 1\circ f_a)(a)\le 1$. Thus we have 
$$(f_a\circ 1\circ f_a)(a)=1.$$

If $A_a=\emptyset$, then $((f_a\circ 1)\circ f_a)(a)=0$ which is 
impossible. Thus we have $A_a\not=\emptyset$. Then\[(({f_a} \circ 1) 
\circ {f_a})(a) = \bigvee\limits_{(x,y) \in {A_a}} {\min \{ ({f_a}}  
\circ 1)(x),{f_a}(y)\} .\]Then there exists $(x,y)\in A_a$ such that 
$(f_a\circ 1)(x)\not=0$ and $f_a(y)\not=0$. Indeed, if $(f_a\circ 
1)(x)=0$ or $f_a(y)=0$ for every $(x,y)\in A_a$, then 
$\min\{(f_a\circ 1)(x),f_a(y)\}=0$ for every $(x,y)\in A_a$, then 
$((f_a\circ 1)\circ f_a)(a)=0$ which is impossible.

Since $f_a(y)\not=0$, we have $y=a$, then $(x,a)\in A_a$. Since 
$(f_a\circ 1)(x)\not=0$, we get $A_x\not=\emptyset$. Since 
$A_x\not=\emptyset$, we have\[({f_a} \circ 1)(x) = 
\bigvee\limits_{(b,c) \in {A_x}} {\min \{ {f_a}} (b),1(c)\}  = 
\bigvee\limits_{(b,c) \in {A_x}}f_a(b).\]If $b\not=a$ for every 
$(b,c)\in A_x$, then $f_a(b)=0$ for every $(b,c)\in A_x$, then 
$(f_a\circ 1)(x)=0$ which is impossible. Hence there exists $(b,c)\in 
A_x$ such that $b=a$. Then $(a,c)\in A_x$. We have $(x,a)\in A_a$ and 
$(a,c)\in A_x$. So we obtain $a\in x\circ a$ and $x\in a\circ c$. 
Then, by Lemma 2, we have
$$a\in x\circ a=\{x\}*\{a\}\subseteq (a\circ c)*\{a\}.$$Since $c\in 
H$ and $a\in (a\circ c)*\{a\}$, the hypersemigroup $H$ is regular. 
$\hfill\Box$ \medskip

\noindent{\bf Definition 6.} An hypersemigroup $H$ is called {\it 
intra-regular} if for every $a\in H$ there exist $x,y\in H$ such 
that$$a\in (x\circ a)*(a\circ y).$$Equivalent Definitions:

(1) $a\in H*\{a\}*\{a\}*H$ for every $a\in H$.

(2) $A\subseteq H*A*A*H$ for every $A\subseteq H$.\medskip

\noindent{\bf Theorem 7.} {\it An hypersemigroup H is intra-regular 
if and only if, for any fuzzy subset f of H, we have$$f\preceq 1\circ 
f\circ f\circ 1.$$}{\bf Proof.} $\Longrightarrow$. Let $f$ be a fuzzy 
subset of $H$ and $a\in H$. Since $H$ is regular, there exist $x,y\in 
H$ such that $a\in (x\circ a)*(a\circ y)$. Then, By Proposition 1(1), 
there exist $u\in x\circ a$ and $v\in a\circ y$ such that $a\in 
u\circ v$. Since $(u,v)\in A_a$, we have\begin{eqnarray*}(1 \circ f 
\circ f \circ 1)(a)&=&\bigvee\limits_{(y,z) \in {A_a}} \min\{(1 \circ 
f)(y),(f \circ 1)(z)\}\\&\ge&\min\{(1 \circ f)(u),(f\circ 1)(v)\}. 
\end{eqnarray*}Since $(x,a)\in A_u$, we have
\[(1 \circ f)(u) = \bigvee\limits_{(y,z) \in {A_u}} {\min \{ 
1(y),f(z)\}  \ge \min \{ 1(x),f(a)\}  = f(a).} \]Since $(a,y)\in 
A_v$, we have\[(f \circ 1)(v) = \bigvee\limits_{(y,z) \in {A_v}} 
{\min \{ f(y),1(z)\}  \ge \min \{ f(a),1(y)\}  = f(a).} \]Hence we 
obtain$$(1\circ f\circ f\circ 1)(a)\ge \min\{f(a),f(a)\}=f(a),$$so $
f\preceq 1\circ f\circ f\circ 1.$\smallskip

\noindent$\Longleftarrow$. Let $a\in H$. Since $f_a$ is a fuzzy 
subset of $H$, by hypothesis, we have $1=f_a(a)\le (1\circ f_a\circ 
1\circ f_a)(a)\le 1$, thus $(1\circ f_a\circ 1\circ f_a)(a)=1$. If 
$A_a=\emptyset$, then $(1\circ f_a\circ 1\circ f_a)(a)=0$, 
impossible, thus $A_a\not=\emptyset$. Then$${\Big(}(1 \circ {f_a}) 
\circ ({f_a} \circ 1){\Big)}(a) = \bigvee\limits_{(x,y) \in {A_a}} 
{\min \{ (1 \circ {f_a}} )(x),({f_a} \circ 1)(y)\}.$$Then there 
exists $(x,y)\in A_a$ such that $(1\circ f_a)(x)\not=0$ and 
$(f_a\circ 1)(y)\not=0$ (otherwise, $(1\circ f_a\circ 1\circ 
f_a)(a)=0$ which is impossible). If $A_x=\emptyset$, then $(1\circ 
f_a)(x)=0$ this is no possible, thus $A_x\not=\emptyset$. Then\[(1 
\circ {f_a}) \circ (x) = \bigvee\limits_{(b,c) \in {A_x}} {\min \{ 
1(b),{f_a}} (c)\}  = \bigvee\limits_{(b,c) \in {A_x}}f_a (c).\]If 
$c\not=a$ for every $(b,c)\in A_x$, then $f_a(c)=0$ for every 
$(b,c)\in A_x$, then $(1\circ f_a)(x)=0$ which is impossible. Then 
there exists $(b,c)\in A_x$ such that $c=a$, so we get $(b,a)\in 
A_x$. If $A_y=\emptyset$, then $(f_a\circ 1)(y)=0$, impossible, thus 
$A_y\not=\emptyset$, and so\[({f_a} \circ 1) \circ (y) = 
\bigvee\limits_{(u,d) \in {A_y}} {\min \{ {f_a}} (u),1(d)\}  = 
\bigvee\limits_{(u,d) \in {A_y}}f_a(u).\]If $u\not=a$ for every 
$(u,d)\in A_y$, then $f_a(u)=0$ for every $(u,d)\in A_y$, and then 
$(f_a\circ 1)(y)=0$ which is no possible. Thus there exists $(u,d)\in 
A_y$ such that $u=a$, then $(a,d)\in A_y$. We have $(x,y)\in A_a$, 
$(b,a)\in A_x$, $(a,d)\in A_y$, that is,$$a\in x\circ y,\; x\in 
b\circ a\;\mbox { and }\;y\in a\circ d.$$Then $a\in x\circ 
y=\{x\}*\{y\}\subseteq (b\circ a)*(a\circ d)$, where $b,d\in H$, so 
the hypersemigroup $H$ is intra-regular. $\hfill\Box$\medskip

\noindent{\bf Definition 8.} An hypersemigroup $H$ is called {\it 
left quasi-regular} if for every $a\in H$ there exist $x,y\in H$ such 
that $$a\in (x\circ a)*(y\circ a).$$Equivalent Definitions:

(1) $a\in H*\{a\}*H*\{a\}$ for every $a\in H$.

(2) $A\subseteq H*A*H*A$ for every $A\subseteq H$.\medskip

\noindent{\bf Theorem 9.} {\it An hypersemigroup H is left 
quasi-regular if and only if, for any fuzzy subset f of H, we 
have$$f\preceq 1\circ f\circ 1\circ f.$$}{\bf Proof.} 
$\Longrightarrow$. Let $a\in H$. Then $f(a)\le (1\circ f\circ 1\circ 
f)(a)$. In fact: Since $H$ is left quasi-regular, there exist $x,y\in 
H$ such that $a\in (x\circ a)*(y\circ a)$. Then there exist $u\in 
x\circ a$ and $v\in y\circ a$ such that $a\in u\circ v$. Since $a\in 
u\circ v$, we have $(u,v)\in A_a$. Since $(u,v)\in A_a$, $A_a$ is a 
nonempty set and we have\begin{eqnarray*}(1 \circ f \circ 1 \circ 
f)(a)&:=&\bigvee\limits_{(y,z) \in {A_a}} \min\{(1 \circ f)(y),(1 
\circ f)(z)\}\\&\ge&\min\{(1 \circ f)(u),(1\circ f)(v)\}. 
\end{eqnarray*}Since $u\in x\circ a$, we have $(x,a)\in A_u$. Then 
$A_u$ is a nonempty set and we have$$(1 \circ 
f)(u):=\bigvee\limits_{(s,t) \in {A_u}} \min\{(1(s), f(t)\ge 
\min\{1(x), f(a)\}=f(a).$$Since $v\in y\circ a$, we have $(y,a)\in 
A_v$. Then $A_v$ is a nonempty set and we have$$(1 \circ 
f)(v):=\bigvee\limits_{(l,k) \in {A_v}} \min\{(1(l), f(k)\ge 
\min\{1(y), f(a)\}=f(a).$$Thus we get$$(1\circ f\circ 1\circ f)(a)\ge 
\min\{f(a),f(a)\}=f(a),$$so $f\preceq 1\circ f\circ 1\circ 
f$.\smallskip

\noindent$\Longleftarrow$. Let $a\in H$. By hypothesis, we have 
$1=f_a(a)\le (1\circ f_a\circ 1\circ f_a)(a)\le 1$, so$$(1\circ 
f_a\circ 1\circ f_a)(a)=1.$$If $A_a=\emptyset$, then $(1\circ 
f_a\circ 1\circ f_a)(a)=0$ which is impossible. Thus we have 
$A_a\not=\emptyset$. Then$${\Big(}(1 \circ {f_a}) \circ (1\circ {f_a} 
){\Big)}(a) = \bigvee\limits_{(x,y) \in {A_a}} {\min \{ (1 \circ 
{f_a}} )(x),(1\circ f_a)(y)\}.$$Then there exists $(x,y)\in A_a$ such 
that $(1\circ f_a)(x)\not=0$ and $(1\circ f_a)(y)\not=0$ (otherwise 
$(1\circ f_a\circ 1\circ f_a)(a)=0$ which is impossible). If 
$A_x=\emptyset$, then $(1\circ f_a)(x)=0$ this is no possible, thus 
$A_x\not=\emptyset$. Then\[(1 \circ {f_a})(x) = \bigvee\limits_{(b,c) 
\in {A_x}} {\min \{ 1(b),{f_a}} (c)\}  = \bigvee\limits_{(b,c) \in 
{A_x}}f_a(c).\]If $c\not=a$ for every $(b,c)\in A_x$, then $f_a(c)=0$ 
for every $(b,c)\in A_x$, then $(1\circ f_a)(x)=0$ which is 
impossible. Then there exists $(b,c)\in A_x$ such that $c=a$. Then we 
have $(b,a)\in A_x$. If $A_y=\emptyset$, then $(1\circ f_a)(y)=0$ 
which is impossible. Thus $A_y\not=\emptyset$. Then$$(1 \circ f_a)(y) 
= \bigvee\limits_{(c,d) \in {A_y}} {\min \{ 1(c),{f_a}} (d)\}  = 
\bigvee\limits_{(c,d) \in {A_y}}f_a(d).$$If $d\not=a$ for every 
$(c,d)\in A_y$, then $(1\circ f_a)(y)=0$ which is impossible. Thus 
there exist $(c,d)\in A_y$ such that $d=a$. Thus we get $(c,a)\in 
A_y$.

We have $(x,y)\in A_a$, $(b,a)\in A_x$, $(c,a)\in A_y$, that is $a\in 
x\circ y$, $x\in b\circ a$, $y\in c\circ a$. Thus we have$$a\in 
x\circ y=\{x\}*\{y\}\subseteq (b\circ a)*(c\circ a),$$ where $b,c\in 
H$, so $H$ is left quasi-regular. $\hfill\Box$\medskip

\noindent{\bf Definition 10.} An hypersemigroup $H$ is called {\it 
right quasi-regular} if for every $a\in H$ there exist $x,y\in H$ 
such that $$a\in (a\circ x)*(a\circ y).$$Equivalent Definitions:

(1) $a\in \{a\}*H*\{a\}*H$ for every $a\in H$.

(2) $A\subseteq A*H*A*H$ for every $A\subseteq H$.\\The right 
analogue of the above theorem also holds, and we have the following 
theorem.\medskip

\noindent{\bf Theorem 11.} {\it An hypersemigroup H is right 
quasi-regular if and only if, for any fuzzy subset f of H, we 
have$$f\preceq f\circ 1\circ f\circ 1.$$}
{\bf Definition 12.} An hypersemigroup $H$ is called {\it semisimple} 
if for every $a\in H$ there exist $x,y,z\in H$ such that $$a\in 
(x\circ a)*(y\circ a)*\{z\}.$$Equivalent Definitions:

(1) $a\in H*\{a\}*H*\{a\}*H$ for every $a\in H$.

(2) $A\subseteq H*A*H*A*H$ for every $A\subseteq H$.\medskip

\noindent{\bf Theorem 13.} {\it An hypersemigroup H is semisimple if 
and only if, for any fuzzy subset f of H, we have$$f\preceq 1\circ 
f\circ 1\circ f\circ 1.$$}{\bf Proof.} $\Longrightarrow$. Let $a\in 
H$. Since $H$ is semisimple, there exist $x,y\in H$ such that $a\in 
(x\circ a)*(y\circ a)*\{z\}={\Big(}(x\circ a)*\{y\}{\Big)}*(a\circ 
z)$. By Proposition 1(1), there exist $u\in (x\circ a)*\{y\}$ and 
$v\in a\circ z$ such that $a\in u\circ v$. Since $u\in (x\circ 
a)*\{y\}$, By Proposition 1(1), there exists $w\in x\circ a$ such 
that $u\in w\circ y$. Thus we have$$a\in u\circ v,\; u\in w\circ y,\; 
w\in x\circ a,\; v\in a\circ z.$$Since $a\in u\circ v$, we have 
$(u,v)\in A_a$. Then $A_a\not=\emptyset$ and\begin{eqnarray*}(1 \circ 
f \circ 1 \circ f\circ 1)(a)&=&\bigvee\limits_{(c,d) \in {A_a}} 
\min\{(1 \circ f\circ 1)(c),(f \circ 1)(d)\}\\&\ge&\min\{(1 \circ 
f\circ 1)(u),(f\circ 1)(v)\}. \end{eqnarray*}Since $u\in w\circ y$, 
we have $(w,y)\in A_u$. Then $A_u\not=\emptyset$ and$$(1 \circ f 
\circ 1)(u)=\bigvee\limits_{(s,t) \in {A_u}} \min\{(1 \circ 
f)(s),1(t)\}\ge(1\circ f)(w).$$Since $w\in x\circ a$, we have 
$(x,a)\in A_w$. Then $A_w\not=\emptyset$ and$$(1 \circ 
f)(w)=\bigvee\limits_{(\xi,\zeta) \in {A_w}} \min\{(1(\xi),f(\zeta)\} 
\ge\min\{1(x),f(a)\}=f(a).$$Since $v\in a\circ z$, we have $(a,z)\in 
A_v$. Then $A_v\not=\emptyset$ and$$(f \circ 1) 
(v)=\bigvee\limits_{(k,h) \in {A_v}} 
\min\{f(k),1(h)\}=\bigvee\limits_{(k,h) \in {A_v}} f(k)\ge 
f(a).$$Thus we have$$(1\circ f\circ 1\circ f\circ 1)(a)\ge 
\min\{f(a),f(a)\}=f(a),$$so $f\preceq 1\circ f\circ 1\circ f\circ 
1$.\smallskip

\noindent$\Longleftarrow$. Let $a\in H$. By hypothesis, we have 
$(1\circ f_a\circ 1\circ f_a\circ 1)(a)=1$. Then $A_a\not=\emptyset$ 
and$$(1\circ f_a\circ 1\circ f_a\circ 1)(a) = \bigvee\limits_{(x,y) 
\in {A_a}}\min \{ (1 \circ {f_a}\circ 1 )(x),(f_a\circ 1)(y)\}.$$
Then there exists $(x,y)\in A_a$ such that $(1\circ f_a\circ 
1)(x)\not=0$ and $(f_a\circ 1)(y)\not=0$.\\If $A_y=\emptyset$, then 
$(f_a\circ 1)(y)=0$ which is impossible. Thus $A_y\not=\emptyset$ and 
$$(f_a\circ 1)(y)=\bigvee\limits_{(b,c)\in {A_y}}\min 
\{f_a(b),1(c)\}=\bigvee\limits_{(b,c)\in A_y}f_a(b).$$If $b\not=a$ 
for every $(b,c)\in A_y$, then $f_a(b)=0$ for every $(b,c)\in A_y$, 
then $(f_a\circ 1)(y)=0$ which is impossible. Thus there exists 
$(b,c)\in A_y$ such that $b=a$, then $(a,c)\in A_y$. If 
$A_x=\emptyset$, then $(1\circ f_a\circ 1)(x)=0$, impossible. Thus 
$A_x\not=\emptyset$ and$$(1\circ f_a\circ 1)(x)= 
\bigvee\limits_{(\rho,\lambda)\in A_x}\min\{1(\rho), (f_a\circ 
1)(\lambda)\}= \bigvee\limits_{(\rho,\lambda)\in A_x}(f_a\circ 
1)(\lambda).$$If $(f_a\circ 1)(\lambda)=0$ for every 
$(\rho,\lambda)\in A_x$, then $(1\circ f_a\circ 1)(x)=0$, impossible. 
Then there exists $(\rho,\lambda)\in A_x$ such that $(f_a\circ 
1)(\lambda)\not=0$. If $A_\lambda=\emptyset$, then $(f_a\circ 
1)(\lambda)=0$, impossible. Thus $A_\lambda\not=\emptyset$ and
\[({f_a} \circ 1)(\lambda)=\bigvee\limits_{(k,h) \in {A_\lambda}} 
{\min \{ {f_a}} (k),1(h)\}=\bigvee\limits_{(k,h) \in {A_\lambda}} 
f_a(k).\]If $a\not=k$ for every $(k,h)\in A_\lambda$, then $(f_a\circ 
1)(\lambda)=0$, impossible. Thus there exists $(k,h)\in A_\lambda$ 
such tha $a=k$, then $(a,h)\in A_\lambda$. We have$$a\in x\circ y,\; 
y\in a\circ c,\; x\in \rho\circ \lambda,\; \lambda\in a\circ h.$$Then 
we have\begin{eqnarray*}a\in 
\{x\}*\{y\}&\subseteq&\{\rho\}*\{\lambda\}*\{a\}*\{c\}\subseteq 
\{\rho\}*\{a\}*\{h\}*\{a\}*\{c\}\\&=&(\rho\circ a)*(h\circ a)*\{c\}, 
\end{eqnarray*} where $\rho,h,c\in H$, so $H$ is semisimple. 
$\hfill\Box$\medskip

\noindent{\bf Note.} The characterization of regular and 
intra-regular hypersemigroups in terms of fuzzy sets given in this 
paper are very useful for further investigation. Exactly as in 
semigroups, using these definitions, many proofs on hypersemigroups 
can be drastically simplified. Let us just give an example to clarify 
what we say, further interesting information concerning this 
structure will be given in a forthcoming paper. We begin with the 
definition of fuzzy right and fuzzy left ideals of hypersemigroups. 
If $H$ is an hypersemigroup, a fuzzy subset of $H$ is called a {\it 
fuzzy right ideal} of $H$ if $f(x\circ y)\ge f(x)$ for every $x,y\in 
H$, in the sense that if $x,y\in H$ and $u\in x\circ y$, then 
$f(u)\ge f(x)$. A fuzzy subset of $H$ is called a {\it fuzzy left 
ideal} of $H$ if $f(x\circ y)\ge f(y)$ for every $x,y\in H$, that is, 
if $x,y\in H$ and $u\in x\circ y$, then $f(u)\ge f(y)$. A fuzzy 
subset $f$ of $H$ is a fuzzy right (resp. fuzzy left) ideal of $H$ if 
and only if $f\circ 1\preceq f$ (resp. $1\circ f\preceq f$) [3]. 
Using the definition of regular hypersemigroups given in the present 
paper, one can immediately see that if $H$ is a regular 
hypersemigroup then, for every fuzzy right ideal $f$ and every fuzzy 
left ideal $g$ of $H$, we have $f\wedge g=f\circ g$. In 
fact,\begin{eqnarray*}f\wedge g&\preceq&(f\wedge g)\circ 1\circ 
(f\wedge g)\preceq (f\circ 1)\circ g\preceq f\circ g\\&\preceq& 
(f\circ 1)\wedge (1\circ g)\preceq f\circ g,\end{eqnarray*}thus 
$f\wedge g=f\circ g$.
{\small\bigskip

\bigskip

\noindent University of Athens, Department of Mathematics\\15784 
Panepistimiopolis, Athens, Greece\\email: nkehayop@math.uoa.gr

\end{document}